\documentclass[a4paper,12pt]{amsart}

\usepackage{amssymb,amsbsy,amsmath,amsfonts,amssymb,amscd}
\usepackage{latexsym}
\usepackage{graphics}
\usepackage{color}
\usepackage{comment}
\input xy
\xyoption{all}

\theoremstyle{plain}
\newtheorem{thm}{Theorem}[section]

\theoremstyle{definition}

\newtheorem{example}[thm]{Example}

\numberwithin{equation}{section}

\newcommand{\sR}{{\mathcal R}}

\newcommand{\PP}{\ensuremath{\mathbb{P}}}

\newcommand{\CC}{\ensuremath{\mathbb{C}}}

\newcommand{\QQ}{\ensuremath{\mathbb{Q}}}

\newcommand{\hol}{\ensuremath{\mathcal{O}}}
\newenvironment{dedication}
        {\begin{quotation}\begin{center}\begin{em}}
        {\par\end{em}\end{center}\end{quotation}}

\newcommand\la{\lambda}

\def\eea{\end{eqnarray*}}
\def\bea{\begin{eqnarray*}}

\newcommand\dual{\mathrel{\raise3pt\hbox{$\underline{\mathrm{\thinspace d
\thinspace}}$}}}
\newcommand\qe{\ifhmode\unskip\nobreak\fi\quad $\Box$}       

\def\BOX{\hfill\lower.5\baselineskip\hbox{$\Box$}}

\newtheorem{theo}{Theorem}[section]

\newtheorem{remarkk}[theo]{Remark}
\newenvironment{rem}{\begin{remarkk}\rm}{\end{remarkk}}

\newenvironment{ex}{\begin{example}\rm}{\end{example}}
\newtheorem{conj}[theo]{Conjecture}

\usepackage{hyperref}
\title [Pluricanonical Maps]{ Pluricanonical Maps and the Fujita Conjecture}

\author{Fabrizio Catanese}
\address {Mathematisches Institut der Universit\"at Bayreuth\\
NW II,  Universit\"atsstr. 30\\
95447 Bayreuth}
\email{fabrizio.catanese@uni-bayreuth.de}
\address{  Korea Institute for Advanced Study, Hoegiro 87, Seoul, 
133--722.}
 
 \date{\today}

\begin{document}

\maketitle

\begin{dedication}
\end{dedication}

\begin{abstract}
We describe examples showing the sharpness of Fujita's conjecture on adjoint bundles also in the general type case, and
use these examples  to formulate  related bold conjectures  on pluricanonical maps of varieties of general type.
\end{abstract}


\section*{Introduction and history of the problem.}

The celebrated Fujita's conjecture on adjoint bundles
\cite{fujita} predicts that, if $H$ is an ample divisor on a smooth projective variety $X$ of dimension $n$, 
then $ D : = K_X + (n+1) H$ is  spanned (that is, $ \hol_X(D)$ is generated by global sections), and 
$K_X + (n+2) H$ is very ample. Fujita observed that, if $X$ is projective space $\PP^n$, and
$H$ is the divisor of a hyperplane,  then $K_X = - (n+1) H$ hence  his conjecture is sharp.

\bigskip

I describe here some  series of  examples that I discovered  some 20 years ago, \footnote{and  which I communicated verbally to several people}
 which show the sharpness of the celebrated Fujita's conjecture on adjoint bundles
also in the case where $X$ is a variety of general type with ample canonical divisor $K_X$ (here we take  $H := K_X$).

\bigskip

The motivation for this short note came from reading the review of the paper \cite{meng}, which shows the  birationaliy of the pentacanonical map
for threefolds  of general type with a Gorenstein  minimal model.

Our examples are  constructed just as an extension  to higher dimensions of some  surfaces classified by Enriques in Chapter VIII of his book
\cite{superficie}, and which, by work of Bombieri \cite{cm}, are characterized exactly as the surfaces 
of general type with the worst possible behaviour
of the pluricanonical maps.

  These examples show the sharpness of some bold question-conjectures which we propose here,
which are easily seen to be some sort of generalization  of   the 
Fujita conjecture \cite{fujita} on adjoint bundles   in the general type case  (part b) below  is due to Meng Chen).

\begin{conj}\label{conj}
Let $X$ be a variety of general type,  which is minimal,  of dimension $n$ and 
with at worst $\QQ$-factorial terminal singularities and such that the canonical divisor $K_X$ is a Cartier divisor:
 then 
 
 a) the pluricanonical map associated to $|m K_X|$ is birational onto its image as soon as $ m \geq n+3$, 
 
 a') the $m$-th pluricanonical map yields an embedding of the canonical model  as soon as $ m \geq n+3$

  b) if $K_X^n = Vol(X) \geq 2$, then the pluricanonical map associated to $|m K_X|$ is birational onto its image as soon as $ m \geq n+2$.
  
  b') if $K_X^n = Vol(X) \geq 2$, then the $m$-th pluricanonical map yields an embedding of the canonical model  as soon as $ m \geq n+2$.
\end{conj}

\begin{rem}

(I) 
Assume that $X$ is smooth of  dimension $n$ and that $K_X$ is an ample divisor. 

Then the above conjecture is implied by the Fujita conjecture on adjoint bundles:
if $H$ is an ample divisor, then $ K_X + (n+2) H$
is very ample.

\medskip

(II) Alternatively, one may formulate the conjecture with the assumption   that $X$ is the canonical model of a variety of
general type, and that $K_X$ is a Cartier divisor.
\end{rem}

  Recall  that Hacon and McKernan \cite{hacmk} proved that for each dimension $n$ there is a constant $c(n)$ such that
for each variety $X$ of general type of dimension $n$ the pluricanonical map associated to $|m K_X|$ is birational onto its image as soon as $ m \geq c(n)$. 

Moreover, we know from \cite{bchm} that the canonical ring 
$$\sR(X) = \oplus _m H^0(X, \hol_X(m K_X))$$
of a variety of general type is finitely generated, hence there exists the canonical model $Y: = Proj (\sR(X))$,
which has canonical singularities, and the volume  $Vol (X)$ is also equal to $K_Y^n$
($Y$ is $\QQ$-Gorenstein).
We also know from ibidem that $X$ admits a minimal model, that is, a model with terminal singularities,
in particular   $\QQ$-Gorenstein. 

In general, see \cite{reid} page 359 for examples in dimension $3$, the question to determine $c(n)$ is complicated by the
fact that the volume need not be an integer: for this reason one should first of  all restrict oneself to  the  case where $Vol(X)$ is an integer.

This is not enough, as observed by Chen Jiang, who  gave the following example.

\begin{ex}
(Chen Jiang) 

Take X to be the Reid-Fletcher hypersurface $X_{46} \subset  \PP(4,5,6,7,23)$, which has volume 
$Vol(X) = \frac{1}{420}$  and  canonical stability index 27. 

Take further a smooth curve $C$ of genus 421. 

Let $ Y: = X \times  C$. 
Then $Y$  is a minimal 4-fold of general type with  canonical volume
$Vol(Y)=8$, an integer. However $7K_Y$  does not yield a  birational map.   More high dimensional examples
can be produced by the same procedure. 
\end{ex}

Our  first series of examples  are varieties of general type which show that both  conjectures are sharp: questions a) a') and 
 the Fujita conjecture (for the latter,  we have mentioned  that it is sharp  for projective space, but  it is useful to produce other
examples, of general type, where it is sharp).

Indeed,  for even dimension $n$, we get smooth varieties $X$ for which $K_X$ is ample and $ K_X + (n+1) H = (n+2) K_X$
is not birational onto its image.

For odd dimension instead  $  (n+2) K_X$ is not birational onto its image, but here $X$ is canonical with terminal singularities
which are not Gorenstein, they are  only 2-Gorenstein.

In particular, we get   threefolds of general type with volume equal to $1$, such that  the $5$-canonical map is not birational.
 The main result of \cite{meng}
is that the $5$-canonical map is  birational if $X$ is Gorenstein.   

 The result shows the subtlety of the higher dimensional geometry, 
where minimal models need not be  Gorenstein:  as observed by Meng Chen, under the Gorenstein assumption, the volume must be even
because of the Riemann-Roch theorem, hence $Vol(X) \geq 2$ and the result confirms part b) of the conjecture.

Our `conjecture' is a theorem obviously for dimension $n=1$, here the smooth curves of genus $2$ show that
the bicanonical map need not be birational, while the tricanonical map is always an embedding.

For surfaces, as already mentioned, the conjecture is again a theorem, due to Bombieri \cite{cm}, and the surfaces
which show that the 4-canonical map need not be  a birational  embedding are exactly the double covers of a quadric cone
branched on the section with a degree 5 surface (these examples are to be found in section 14 of Chapter VIII, pages 303-305
of \cite{superficie}) while the other surfaces for which the 3-canonical map need not be  a birational  embedding  are
the double covers of the plane branched on a curve of degree 8 ( see section 17, Chapter VIII,  pages 311-312 of \cite{superficie}).

 Our second series of examples are smooth varieties of general type $X$ with $K_X$ ample, and show that conjectures 
b), b') are sharp. They also show that for Fujita's conjecture one may also ask whether $ K_X + (n+1)H$
is very ample once we know that $ H^n \geq 2$.

 In higher dimension, the  general case of varieties of general type is quite complicated:
already  for  threefolds of general type whose minimal  model is not Gorenstein the situation is intricate:
improving on  \cite{chen2}, \cite{mengnew}   proves birationality for $ m \geq 57$, else one needs some other hypotheses,
for instance  the hypothesis
of positive  irregularity \cite{meng2} \cite{chen3}.

Also concerning  the Fujita conjecture there are until now no counterexamples, and  the best result
was obtained by Angehrn and Siu \cite{siu}, showing that $K_X + m H$ is very ample for
$ m \geq \frac{1}{2} (n^2 + n +2)$.

\section{A series of hypersurfaces in weighted projective space}

For simplicity we stick to algebraic varieties defined over $\CC$, since the above conjectures are meant
for geometry over the complex numbers: but the examples make perfect sense over any algebraically
closed field $k$  of  $char (k) \neq 2$ ($char(k) = 2$   requires some minor modification).

\begin{example}
For each $n$, consider the weighted projective space $\PP (1^n, 2, n+3)$, with coordinates 
$$(x,y,z) : = (x_1, \dots, x_n, y, z),$$
and choose $X$ to be a general hypersurface of degree $ 2 (n+3)$,
$$ X : = \{ (x,y,z) | z^2 = F (x,y)\}.$$
\begin{enumerate}
\item
$X$ is a double cover of $\PP (1^n, 2)$, which is the projective cone over the quadratic Veronese embedding
of $\PP^{n-1}$, which therefore is smooth, except at the point $x=0$, where we have the `Kummer' singularity
$\CC^n / \pm 1.$ 
The singularity is Gorenstein if and only if $n$ is even; it is terminal for $n \geq 2$, since the minimal resolution
is the quotient of the blow up of $\CC^n$ at the origin, hence it carries an exceptional divisor 
$ E \cong \PP^{n-1}$ with normal bundle $ = \hol_{\PP^{n-1}}(-2)$. Therefore the canonical divisor of the resolution
is the pull back of the canonical divisor of $X$ plus $ - \frac{1}{2} (n-2) E$. 
\item
The double covering is branched on the hypersurface $ F (x,y)=0$, which for general choice of $F$
is smooth and does not contain the singular point of $\PP (1^n, 2)$. Hence $X$ is smooth except possibly for the
points lying above the point $x=0$.
\item
$X$ is generally smooth for $n$ even, but it has two Kummer singularities if $n$ is odd.
In fact, over the point $x=0$, we get the points defined by $z^2 = \la y^{n+3}$, where $\la$ is the coefficient of  
$y^{n+3}$ in $F$; there seems to be two distinct solutions for $\la \neq 0$, $y=1, z= \pm 1$,
but, for $n$ even, $n+3$ is odd, hence the diagonal action of $-1$ sends $y \mapsto y, z \mapsto -z$.
Hence, we get two Kummer singularities if $n$ is odd, but for $n$ even we get a double covering ramified exactly at the singular point,
hence we get a smooth point.
\item
Adjunction (see \cite{dolgachev})  shows that $\omega_X = \hol_X(1)$, hence for $X$ even we have smooth variety with 
$K_X$ ample; for $n$ odd we have that $X$ is terminal and $2$-Gorenstein, as seen above.
\item
$H^0( \hol_X(m))$ yields a birational map if and only if $ m \geq n+3$, which is indeed an embedding.
\end{enumerate}
\end{example} 

The next series extends to higher dimension the other surfaces considered by Enriques, which require $ m \geq 4$:

\begin{ex}
For each $n$, consider the weighted projective space $\PP (1^{n+1}, n+2)$, with coordinates 
$$(x,z) : = (x_1, \dots, x_n, x_{n+1}, z),$$
and choose $X$ to be a general hypersurface of degree $ 2 (n+2)$,
$$ X : = \{ (x,z) | z^2 = F (x)\}.$$
\begin{enumerate}
\item
$X$ is a double cover of $\PP^n$,  branched on a hypersurface of degree $2n + 4$,
hence 
$X$ is generally smooth.
\item
Adjunction (see \cite{dolgachev})  shows that $\omega_X = \hol_X(1)$, hence  we have smooth variety with 
$K_X$ ample; \item
$H^0( \hol_X(m))$ yields a birational map if and only if $ m \geq n+2$, which is indeed an embedding.
\end{enumerate}
\end{ex}

\section{Post Scriptum}

The Fujita conjecture is true for reduced curves, but not for non-reduced curves which are not numerically connected: 
take an elliptic 
curve $E$ in a smooth algebraic surface $S$
with selfintersection $E^2 = -d, d \geq 3$,  set $C: =2D$, and take a divisor $H$ with $H \cdot E=1$
($S$ may just be the blow up of $E \times \PP^1$).

Then on $C$ we have $ 3 H + K_C = 3 H + K_S + C = 3 H + K_S + 2 E$ and its intersection with $E$ equals $3 + E^2 = 3 - d$,
hence this linear system is not even ample as soon as $ d \geq 3$. 

\bigskip

{\bf Added in proof:} Burt Totaro mentioned to me that similar examples are to be found in 
Theorem 3 of \cite{bpt}, which contains also  other examples 
where the minimal model does not have a Cartier canonical divisor $K_X$.

Under the condition that the minimal models have a canonical divisor  $K_X$ which is not Cartier,
then there are recent more extreme and striking examples by 
Esser, Wang, and Totaro; these examples are such that there are  about $2^{2^{n/2}}$
vanishing plurigenera in dimension $n$, see  Theorem 1.1 of \cite{etw} and ensuing discussion.

\bigskip

{\bf Acknowledgement:} I would like to thank Gavin Brown for trying to answer many years ago my queries about the search for 
other hypersurfaces, respectively complete
intersections, in weighted projective spaces,  which would show the sharpness of Fujita's conjecture or even disprove it.

Many thanks to Meng Chen for very useful comments on a first draft of this note, and especially for contributing part b)
of Conjecture \ref{conj}. Thanks to Chen Jiang for the permission to include his example in this note.
Thanks to the referee and to Burt Totaro for  interesting comments on the first version of this article.

\end{document}